\documentclass[10pt]{amsart}

\usepackage{amssymb}
\usepackage{enumerate}

\numberwithin{equation}{section}

\newtheorem{theorem}{Theorem}
\newtheorem{lemma}{Lemma}
\newtheorem{proposition}{Proposition}
\newtheorem{corollary}{Corollary}
\newtheorem{definition}{Definition}

\begin{document}

\title[Generalized Benjamin-Ono equation]{Sharp well-posedness results for the generalized Benjamin-Ono equation with high nonlinearity}
\author{St\'ephane Vento}
\address{Universit\'e de Marne-La-Vall\'ee, Laboratoire d'Analyse et
de Math\'ematiques Appliqu\'ees,\\ 5, bd. Descartes, Cit\'e
Descartes, Champs-Sur-Marne,\\ 77454 Marne-La-Vall\'ee Cedex 2,
France.}

\email{stephane.vento@univ-mlv.fr}
\date{}

\def\sgn{\textrm{sgn}}
\def\R{\mathbb{R}}
\def\Z{\mathbb{Z}}
\def\S{\mathcal{S}}
\def\H{\mathcal{H}}
\def\eps{\varepsilon}
\def\pv{\textrm{pv}}
\def\supp{\textrm{supp}}

\begin{abstract}
We establish the local well-posedness of the generalized
Benjamin-Ono equation $\partial_tu+\mathcal{H}\partial_x^2u\pm
u^k\partial_xu=0$  in $H^s(\R)$, $s>1/2-1/k$ for $k\geq 12$ and
without smallness assumption on the initial data. The condition
$s>1/2-1/k$ is known to be sharp since the solution map
$u_0\mapsto u$ is not of class $\mathcal{C}^{k+1}$ on $H^s(\R)$
for $s<1/2-1/k$. On the other hand, in the particular case of the
cubic Benjamin-Ono equation, we prove the ill-posedness in
$H^s(\R)$, $s<1/3$.
\end{abstract}

\maketitle

\section{Introduction and statement of the results}
\subsection{Introduction}

Our purpose in this paper is to study the initial value problem
for the generalized Benjamin-Ono equation
\begin{equation}\label{gBO}\tag{gBO}\left\{\begin{array}{ll}\partial_tu+\mathcal{H}\partial_x^2u\pm u^k\partial_xu=0,
\quad x,t\in\mathbb{R},\\u(x,t=0)=u_0(x),\quad
x\in\mathbb{R},\end{array}\right.\end{equation} where
$k\in\mathbb{N}\setminus\{0\}$, $\mathcal{H}$ is the Hilbert
transform defined by
\[\mathcal{H}f(x)=\frac{1}{\pi}\pv\Big(\frac{1}{x}\ast
u\Big)(x)=\mathcal{F}^{-1}\big(-i\ \sgn(\xi)\hat{f}(\xi)\big)(x)\]
and with initial data $u_0$ belonging to the Sobolev space
$H^s(\R)=(1-\partial^2_x)^{-s/2}L^2(\R)$. \vskip 0.5cm The case
$k=1$ was deduced by  T.B. Benjamin \cite{1967JFM....29..559B} and
later by H. Ono \cite{MR0398275} as a model in internal wave
theory. The Cauchy problem for the Benjamin-Ono equation has been
extensively studied. It has been proved in \cite{MR533234} that
(BO) is globally well-posed (i.e. global existence, uniqueness and
persistence of regularity of the solution) in $H^s(\R)$ for $s\geq
3$, and then for $s\geq 3/2$ in \cite{MR1097916} and
\cite{MR847994}. Recently, T. Tao \cite{MR2052470} proved the
well-posedness of this equation for $s\geq 1$ by using a gauge
transformation. More recently, combining a gauge transformation
with a Bourgain's method, A.D. Ionescu and C.E. Kenig
\cite{math.AP/0508632} shown that one could go down to $L^2(\R)$,
and this seems to be, in some sense, optimal. It is worth noticing
that all these results have been obtained by compactness methods.
On the other hand, L. Molinet, J.-C. Saut and N. Tzvetkov
\cite{MR1885293}  proved that, for all $s\in\R$, the flow map
$u_0\mapsto u$ is not of class $\mathcal{C}^2$ from $H^s(\R)$ to
$H^s(\R)$. Furthermore, building suitable families of approximate
solutions,  H. Koch and N. Tzvetkov   proved in \cite{MR2172940}
that the flow map is not even uniformly continuous on bounded sets
of $H^s(\R)$, $s>0$. As an important consequence of this, since a
Picard iteration scheme would imply smooth dependance upon the
initial data, one see that such a scheme cannot be used to get
solutions in any space continuously embedded in
$\mathcal{C}([0,T],H^s(\R))$.

\vskip 0.5cm

For higher nonlinearities, that is for $k\geq 2$, the picture is a
little bit different. It turns out that one can get local
well-posedness results throught a Picard iteration scheme but for
small initial data only. This seems mainly due to the fact that
the smoothing properties of the linear group $V(\cdot)$ associated
to the linear (BO) equation is just sufficient to recover the lost
derivative in the nonlinear term, but does not allow to get the
required contraction factors. On the other hand, for large initial
data, one can prove local well-posedness by compactness methods
together with a gauge transformation. Unfortunately, this usually
requires more smoothness on the initial data. We summurize now the
known results about the Cauchy problem for (gBO) equations when
$k\geq 2$.

\vskip 0.3cm

In the case of the modified Benjamin-Ono equation ($k=2$), C.E.
Kenig and H. Takaoka \cite{math.AP/0509573} have recently obtained
the global well-posedness in the energy space $H^{1/2}(\R)$. This
have been proved thanks to a localized gauge transformation
combined with a $L^2_{xT}$ estimate of the solution. This result
is known to be sharp since the solution map $u_0\mapsto u$ is not
$\mathcal{C}^3$ in $H^s(\R)$, $s<1/2$ (see \cite{MR2038121}).
\vskip 0.3cm

For (\ref{gBO}) with cubic nonlinearity ($k=3$), the local
well-posedness is known  in $H^s(\R)$, $s>1/3$ for small initial
data \cite{MR2038121} but only in $H^s(\mathbb{R})$, $s>3/4$, for
large initial data. Moreover, the ill-posedness has been proved in
$H^s(\mathbb{R})$, $s<1/6$ \cite{MR2038121}. In this paper, we
show the ill-posedness of the cubic Benjamin-Ono equation in
$H^s(\mathbb{R})$, $s<1/3$, which turns out to be optimal
according to the above results.

\vskip 0.3cm

When $k\geq 4$, by a scaling argument, one can guess the best
Sobolev space in which the Cauchy problem is locally well-posed,
that is, the critical indice $s_c$ such that (\ref{gBO}) is
well-posed in $H^{s}(\R)$ for $s>s_c$ and ill-posed for $s<s_c$ .
Recall that if $u(x,t)$ is a solution of the equation then
$u_\lambda(x,t)=\lambda^{1/k}u(\lambda x,\lambda^2t)$
($\lambda>0$) solves (\ref{gBO}) with initial data
$u_\lambda(x,0)$ and moreover
\[\|u_\lambda(\cdot,0)\|_{\dot{H}^s}=\lambda^{s+\frac{1}{k}-\frac{1}{2}}\|u(\cdot,0)\|_{\dot{H}^s}.\]
Hence the $\dot{H}^s(\R)$ norm is invariant if and only if
$s=s_k=1/2-1/k$ and one can conjecture that $s_c=s_k$.

In the case of small initial data, this limit have been reached by
L. Molinet and F. Ribaud \cite{MR2038121}. This result is almost
sharp in the sense that the flow map $u_0\mapsto u$ is not of
class $\mathcal{C}^{k+1}$ from $H^s(\R)$ to
$\mathcal{C}([0,T],H^s(\R))$ at the origin when $s<s_k$,
\cite{MR2101982}. This lack of regularity is also described by
H.A. Biagioni and F. Linares in \cite{MR1837253} where they
established, using solitary waves, that the flow map is not
uniformly continuous in $\dot{H}^{s_k}(\R)$, $k\geq 2$.

For large initial data, the local well-posedness of (gBO) is only
known in $H^s(\mathbb{R})$, $s\geq 1/2$, whatever the value of
$k$. This have been proved in \cite{MR2101982} by using the gauge
transformation
\begin{equation}\label{GT}u\stackrel{\mathcal{G}}{\longmapsto} P_+(e^{-i\int_{-\infty}^xu^k}u),\end{equation} together with compactness methods.
 Note also that very recently, in the particular case $k=4$, N. Burq and F. Planchon
\cite{math.AP/0409379} derived  the local well-posedness of (gBO)
in the homogeneous space $\dot{H}^{1/4}(\R)$.

In this paper, our aim is to improve the results obtained in
\cite{MR2101982} for large initial data. We show that for all
$k\geq 12$, (\ref{gBO}) is locally well-posed in $H^s(\R)$,
$s>s_k$.
 Our proofs follow those of
\cite{MR2101982} : we perform the gauge transformation
$w=\mathcal{G}(u)$ of a smooth solution $u$ of (gBO) and derive
suitable estimates for $w$. The main interest of this
transformation is to obtain an equation satisfied by $w$ where the
nonlinearity $u^ku_x$ is replaced by terms of the form
$P_+(u^kP_-u_x)$ in which one can share derivatives on $u$ with
derivatives on $u^k$. Working in the surcritical case, this allows
to get a contraction factor $T^\nu$ in our estimates. It is worth
noticing that $\nu=\nu(s)$ verifies $\lim_{s\rightarrow
s_k}\nu(s)=0$, and this explains why our method fails in the
critical case $s=s_k$. On the other hand, the restriction $k\geq
12$ appears when we estimate the integral term
$$P_+\Big(e^{-i\int_{-\infty}^x u^k}u\int_{-\infty}^xu^{k-2}\H
u_{xx}\Big)$$ (see section 3.2). This term doesn't seem to have a
"good structure" since the bad interaction
$$Q_ju\int_{-\infty}^x(P_j u)^{k-2}\H P_ju_{xx}$$ forbids the share
of the antiderivative $\int_{-\infty}^x$ with other derivatives.

\subsection{Main results}
Our main results read as follows.
\begin{theorem}\label{mainresult} Let $k\geq 12$ and $u_0\in H^s(\R)$ with $s>1/2-1/k$. Then there
exist $T=T(s,k,\|u_0\|_{H^s})>0$ and a  unique solution
$u\in\mathcal{C}([0,T];H^s(\R))$ of (\ref{gBO}) such that
\begin{eqnarray}&&\|D^{s+1/2}_xu\|_{L^\infty_xL^2_T}<\infty,\\
&&\|D^{s-1/4}_xu\|_{L^4_xL^\infty_T}<\infty,\\
&&\|P_0u\|_{L^2_xL^\infty_T}<\infty.\end{eqnarray} Moreover, the
flow map $u_0\mapsto u$ is Lipschitz on every bounded set of
$H^s(\R)$.
\end{theorem}

As mentioned previously, these results are in some sense almost
sharp. However, the critical case $s=s_k$ remains open. We will
only consider the most difficult case, that is the lowest values
for $s$. More precisely we will prove Theorem \ref{mainresult} for
$s_k<s<1/2$.

\vskip 0.5cm

In the case $k=3$, we have the following ill-posedness result.

\begin{theorem}\label{illth} Let $k=3$ and $s<1/3$. There does not exist $T>0$ such that
the Cauchy problem (\ref{gBO}) admits an unique local solution
defined on the interval $[0,T]$ and such that the flow map
$u_0\mapsto u$ is of class $\mathcal{C}^4$ in a neighborhood of
the origin from $H^s(\mathbb{R})$ to $H^s(\mathbb{R})$.
\end{theorem}

This result implies that we cannot solve (\ref{gBO}) with $k=3$ in
$H^s(\R)$, $s<1/3$ by a contraction method on the Duhamel
formulation. Recall that for small initial data \cite{MR2038121},
we have local well-posedness in $H^s(\R)$ for $s>1/3$. In view of
this, we can conjecture that (\ref{gBO}) is locally well-posed in
$H^s(\R)$, $s>1/3$. \vskip 0.3cm The remainder of this paper is
organized as follows. In section 2, we first derive some linear
estimates on the free evolution operator associated to (\ref{gBO})
and we define our resolution space. Then we give some technical
lemmas which will be used for nonlinear estimates. In section 3 we
introduce the gauge transformation and derive the needed nonlinear
estimates. The section 4 is devoted to the proof of Theorem
\ref{mainresult}. Finally we prove our ill-posedness result in the
Appendix.

The author is grateful to Francis Ribaud for several useful
comments on the subject.

\subsection{Notations}
For two positive numbers $x$, $y$, we write $x\lesssim y$ to mean
that there exists a $C>0$ which does not depend on $x$ and $y$,
and such that $x\leq Cy$. In the sequel, this constant may depend
on $s$ and $k$. We also use $\nu=\nu(s,k)$ to denote a positive
power of $T$ which may differ at each occurrence.

Our resolution space is constructed thanks to the space-time
Lebesgue spaces $L^p_xL^q_T$ and $L^q_TL^p_x$ endowed for $T>0$
and $1\leq p,q\leq\infty$ with the norm
$$\|f\|_{L^p_xL^q_T}=\big\|\|f\|_{L^q_T([0;T])}\big\|_{L^p_x(\R)}\ \textrm{ and
}\
\|f\|_{L^q_TL^p_x}=\big\|\|f\|_{L^p_x(\R)}\big\|_{L^q_T([0;T])}.$$
When $p=q$ we simplify the notation by writing $L^p_{xT}$.

The well-known operators $\mathcal{F}$ (or $\hat{\cdot}$) and
$\mathcal{F}^{-1}$ (or $\check{\cdot}$) are the Fourier operators
defined by $\hat{f}(\xi)=\int_\R e^{-ix\xi}f(x)dx$. The
pseudo-differential operator $D^\alpha_x$ is defined by its
Fourier symbol $|\xi|^\alpha$. Let $P_+$ and $P_-$ be the Fourier
projections to $[0,+\infty[$ and $]-\infty,0]$. Thus one has
$$i\mathcal{H}=P_+-P_-.$$
Let $\eta\in\mathcal{C}_0^\infty(\R)$, $\eta\geq 0$, $\supp\
\eta\subset\{1/2\leq |\xi|\leq 2\}$ with
$\sum_{-\infty}^\infty\eta(2^{-k}\xi)=1$ for $\xi\neq 0$. We set
$p(\xi)=\sum_{j\leq -3}\eta(2^{-j}\xi)$ and consider, for all
$k\in\Z$, the operators $Q_k$ and $P_k$ respectively defined by
\[Q_k(f)=\mathcal{F}^{-1}(\eta(2^{-k}\xi)\hat{f}(\xi))\ \textrm{ and
}\ P_k(f)=\mathcal{F}^{-1}(p(2^{-k}\xi)\hat{f}(\xi)).\] Therefore
we have the standard Littlewood-Paley decomposition
\begin{equation}\label{DLP}f=\sum_{j\in\Z}Q_j(f)=P_0(f)+\sum_{j\geq
-2}Q_j(f)=P_0(f)+\Tilde{P}(f).\end{equation} We also need the
operators \[P_{\leq k}f=\sum_{j\leq k}Q_jf,\quad P_{\geq
k}f=\sum_{j\geq k}Q_jf.\] We finally introduce the operators
$\tilde{P}_+=P_+\Tilde{P}$ and $\tilde{P}_-=P_-\Tilde{P}$ in order
to obtain the smooth decomposition
\begin{equation}\label{smoothdec}f=\tilde{P}_-(f)+P_0(f)+\tilde{P}_+(f).\end{equation}

\section{Linear estimates and technical lemmas}
\subsection{Linear estimates and resolution space}

Recall that (\ref{gBO}) is equivalent to its integral formulation
\begin{equation}\label{duhamel}u(t)=V(t)u_0\mp\frac{1}{k+1}\int_0^tV(t-\tau)\partial_x(u^{k+1})(\tau)d\tau,\end{equation}
where $V(t)=\mathcal{F}^{-1}e^{it\xi|\xi|}\mathcal{F}$ is the
generator of the free evolution. Let us now gather the well-known
estimates on the group $V(\cdot)$ in the following lemma.

\begin{lemma}\label{estlin} Let $\varphi\in\mathcal{S}(\R)$, then
\begin{eqnarray}\label{est0}\|V(t)\varphi\|_{L^\infty_TL^2_x} &\lesssim
&\|\varphi\|_{L^2},\\
\label{est1}\|D^{1/2}_xV(t)\varphi\|_{L^\infty_x L^2_T} &\lesssim
&
\|\varphi\|_{L^2},\\
\label{est2}\|D^{-1/4}_xV(t)\varphi\|_{L^4_x L^\infty_T} &\lesssim
& \|\varphi\|_{L^2}.\end{eqnarray} Moreover, for $0<T<1$, we have
\begin{eqnarray}
\label{est4}\|P_0V(t)\varphi\|_{L^2_xL^\infty_T} &\lesssim &
\|P_0\varphi\|_{L^2}.\end{eqnarray}
\end{lemma}
The estimate (\ref{est0}) is straightforward whereas the proof of
the Kato smoothing effect (\ref{est1}) and the maximal in time
inequality (\ref{est2}) can be found in \cite{MR1101221}. Estimate
(\ref{est4}) has been proved in \cite{MR1086966}.

These estimates motivate the definition of our resolution space.

\begin{definition} For $s_k<s<1/2$, we define the
space $X^s_T=\{u\in\mathcal{S}'(\R^2), \|u\|_{X^s_T}<\infty\}$
where $0<T<1$ and
\begin{equation}\label{space}\|u\|_{X^s_T}=\|u\|_{L^\infty_TH^s_x}+\|D^{s+1/2}_xu\|_{L^\infty_xL^2_T}+\|D^{s-1/4}_xu\|_{L^4_xL^\infty_T}
+\|P_0u\|_{L^2_xL^\infty_T}.\end{equation}
\end{definition}

Thus lemma \ref{estlin} implies immediately that for all
$\varphi\in\mathcal{S}(\R)$ and $0<T<1$,
\begin{equation}\label{estsg}\|V(t)\varphi\|_{X^s_T}\lesssim
\|\varphi\|_{H^s}.\end{equation}

We now give some families of norms which are controlled by the
$X^s_T$ norm. This will be usefull to derive some nonlinear
estimates in the sequel.

\begin{definition} A triplet $(\alpha,p,q)\in\R\times[2,\infty]^2$
is said to be
 1-admissible if
$(\alpha,p,q)=(1/2,\infty,2)$ or \begin{equation}\label{1ad} 4\leq
p<\infty,\quad 2< q\leq\infty,\quad
\frac{2}{p}+\frac{1}{q}\leq\frac{1}{2},\quad
\alpha=\frac{1}{p}+\frac{2}{q}-\frac{1}{2}.
\end{equation}
\end{definition}

\begin{proposition}\label{admis} If $(\alpha-s,p,q)$ is 1-admissible, then for all $u$ in
$X^s_T$,
\begin{equation}\label{in1ad}\|D^\alpha_xu\|_{L^p_xL^q_T}\lesssim
\|u\|_{X^s_T}.\end{equation}
\end{proposition}
Proof : The inequality
\begin{equation}\label{triest}\|D^{s+1/2}_xu\|_{L^\infty_xL^2_T}\lesssim\|u\|_{X^s_T}\end{equation}
yields the result when $(\alpha,p,q)=(1/2,\infty,2)$. Assume now
$(\alpha,p,q)\neq (1/2,\infty,2)$. Let $r\in[4;p]$. Then according
to Sobolev embedding theorem,
\[\|D^{s+1/r-1/2}_xu\|_{L^r_xL^\infty_T}\lesssim
\|D^{s-1/4}_xu\|_{L^4_xL^\infty_T}\lesssim \|u\|_{X^s_T}.\] By
interpolation with (\ref{triest}) we get for all $0\leq\theta\leq
1$
\[\|D^{s+\frac{1}{2}-(1-\frac{1}{r})\theta}_xu\|_{L^{r/\theta}_xL^{2/(1-\theta)}_T}\lesssim\|u\|_{X^s_T}.\]
We deduce (\ref{in1ad}) by taking $\theta=r/p$ since the
assumption $r\geq 4$ is equivalent to $\frac{2}{p}+\frac{1}{q}\leq
\frac{1}{2}$. $\Box$

We list now all the norms needed for the nonlinear estimates.

\begin{corollary}\label{bound} For $u\in X^s_T$, the following quantities are bounded by
$\|u\|_{X^s_T}$. $$\begin{array}{ll}
N_1=\|u\|_{L^p_xL^\infty_T},\quad 4\leq
p\leq(\frac{1}{2}-s)^{-1},&
N_2=T^{-\nu}\|u\|_{L^{3k}_{xT}},\\
 N_3=T^{-\nu}\|u\|_{L^{k/(1-s)}_xL^{2k/s}_T}, & N_4=
T^{-\nu}\|u\|_{L^{k(\frac{1}{3}+s)^{-1}}_xL^{k(\frac{1}{3}-\frac{s}{2})^{-1}}_T},\\
N_5=T^{-\nu}\|u\|_{L^{3k/4s}_xL^{k(\frac{1}{2}-\frac{2s}{3})^{-1}}_T},
&
 N_6=T^{-\nu}\|u\|_{L^{k(1-\frac{s}{3})^{-1}}_xL^{6k/s}_T}, \\
 N_7=T^{-\nu}\|u\|_{L^{k+1}_xL^{2k(k+1)}_T}, &
N_{8}=T^{-\nu}\|u\|_{L^{(k-1)(\frac{5}{6}-\frac{s}{3})^{-1}}_xL^{(k-1)(\frac{2s}{3}-\frac{1}{6})^{-1}}_T},\\
N_{9}=\|D^{1-2s+6\eps}_xu\|_{L^{(\frac{3}{2}-3s)^{-1}}_xL^{1/3\eps}_T},
& N_{10}=\|D^s_xu\|_{L^6_{xT}},
\\
N_{11}=\|D^{s+1/2-3\eps}_xu\|_{L^{1/\eps}_xL^{(\frac{1}{2}-2\eps)^{-1}}_T},
&
N_{12}=\|D^{1/2}_xu\|_{L^{3/s}_xL^{(\frac{1}{2}-\frac{2s}{3})^{-1}}_T},
\end{array}$$
where $\eps,\nu>0$ are small enough.
\end{corollary}
Proof : \begin{enumerate} \item[(i)] Let $4\leq p\leq
(\frac{1}{2}-s)^{-1}$. By separating low and high frequencies,
\[\|u\|_{L^p_xL^\infty_T}\lesssim \|P_0u\|_{L^2_xL^\infty_T}+\|\Tilde{P}D^{s+1/p-1/2}_xu\|_{L^p_xL^\infty_T}\lesssim\|u\|_{X^s_T}.\] Here
we used that $\Tilde{P}$ is continuous on $L^p_xL^q_T$, $1\leq
p,q\leq\infty$, and the 1-admissibility of $(1/p-1/2,p,\infty)$.
\item[(ii)-(vii)] We evaluate the norm of the form
$N=\|u\|_{L^p_xL^q_T}$ with $p>2$ and $q<\infty$. Fix $\delta>0$
small enough so that $\alpha=s-s_k-2\delta>0$ and
$\frac{1}{q}-\delta>0$. Then using the previous decomposition,
Bernstein and H\"older inequalities, we get \[N\lesssim
T^\nu\|P_0u\|_{L^2_xL^\infty_T}+T^\nu\|\tilde{P}D^{\alpha}_xu\|_{L^p_xL^{(\frac{1}{q}-\delta)^{-1}}_T}.\]
One complete the proof by noticing that the triplet
$(\alpha-s,p,(\frac{1}{q}-\delta)^{-1})$ is 1-admissible.
\item[(viii)] Following the same idea, we write \[N_{8}\lesssim
T^\nu\|P_0u\|_{L^2_xL^\infty_T}+
T^\nu\|\tilde{P}D^{\frac{k}{k-1}(s-s_k-2\frac{\delta}{k})}_xu\|_{L^{(k-1)(\frac{5}{6}-\frac{s}{3})^{-1}}_x
L^{(k-1)(\frac{2s}{3}-\frac{1}{6}-\delta)^{-1}}_T}\] for an
appropriate $\delta>0$. Once again,
$(\frac{k}{k-1}(s-s_k-2\frac{\delta}{k})-s,(k-1)(\frac{5}{6}-\frac{s}{3})^{-1},(k-1)(\frac{2s}{3}-\frac{1}{6}-\delta)^{-1})$
is 1-admissible. \item[(ix)-(xii)] Note finally that the triplets
$(1-3s+6\eps,(\frac{3}{2}-3s)^{-1},1/3\eps)$, $(0,6,6)$,
$(1/2-3\eps,1/\eps,(\frac{1}{2}-2\eps)^{-1})$ and
$(1/2-s,3/s,(\frac{1}{2}-\frac{2s}{3})^{-1})$ are 1-admissible.
$\Box$
\end{enumerate}

We now turn to the non-homogenous estimates. Let us first recall
the following result found in \cite{MR2101982}.

\begin{lemma}\label{estsnonhom} Let $(\alpha_1,\alpha_2)\in\R^2$,
$(\nu_1,\nu_2)\in \R_+^2$, and $1\leq p_1,q_1,p_2,q_2\leq \infty$
such that for all $\varphi\in\mathcal{S}(\R)$,
\[\|D^{\alpha_1}_xV(t)\varphi\|_{L^{p_1}_xL^{q_1}_T}\lesssim
T^{\nu_1}\|\varphi\|_{L^2},\]
\[\|D^{\alpha_2}_xV(t)\varphi\|_{L^{p_2}_xL^{q_2}_T}\lesssim
T^{\nu_2}\|\varphi\|_{L^2}.\] Then for all
$f\in\mathcal{S}(\R^2)$,
\begin{equation}\label{estnonhom}\Big\|D^{\alpha_2}_x\int_0^tV(t-\tau)f(\tau)d\tau\Big\|_{L^\infty_TL^2_x}\lesssim
T^{\nu_2}\|f\|_{L^{\bar{p}_2}_xL^{\bar{q}_2}_T},\end{equation}
\begin{equation}\label{estnonhom2}\Big\|D^{\alpha_1+\alpha_2}_x\int_0^tV(t-\tau)f(\tau)d\tau\Big\|_{L^{p_1}_xL^{q_1}_T}
\lesssim
T^{\nu_1+\nu_2}\|f\|_{L^{\bar{p}_2}_xL^{\bar{q}_2}_T}\end{equation}
provided $\min(p_1,q_1)>\max(\bar{p}_2,\bar{q}_2)$ or
($q_1=\infty$ and $\bar{p}_2,\bar{q}_2<\infty$), where $\bar{p}_2$
and $\bar{q_2}$ are defined by $1/\bar{p}_2=1-1/p_2$ and
$1/\bar{q}_2=1-1/q_2$.
\end{lemma}

Using lemma \ref{estsnonhom} we infer the following result.

\begin{lemma}\label{lemestnonhom} For all $f\in\mathcal{S}(\R^2)$, the quantity
$\displaystyle\Big\|\int_0^tV(t-\tau)f(\tau)d\tau\Big\|_{X^s_T}$
can be estimated by
\begin{equation}\label{estnonhoms}\|f\|_{L^{(\frac{5}{6}+\frac{s}{3})^{-1}}_xL^{(\frac{5}{6}-\frac{2s}{3})^{-1}}_T},\quad
\|D^s_xf\|_{L^{6/5}_{xT}},\quad
\|D^{s-1/2}_xf\|_{L^1_xL^2_T},\quad
\|D^{s+1/4}_xf\|_{L^{4/3}_xL^1_T}.\end{equation} Moreover,
\begin{equation}\label{estnonhoms2}\Big\|D^{s+1/2}_x\int_0^tV(t-\tau)f(\tau)d\tau\Big\|_{L^\infty_xL^2_T}\lesssim\|D^s_xf\|_{L^1_TL^2_x}.\end{equation}
\end{lemma}

Proof : (\ref{estnonhoms}) follows from
(\ref{estnonhom})-(\ref{estnonhom2}) since the triplets
$(s,(\frac{1}{6}-\frac{s}{3})^{-1},(\frac{1}{6}+\frac{2s}{3})^{-1})$,
$(0,6,6)$, $(1/2,\infty,2)$ and $(-1/4,4,\infty)$ are
1-admissible. Inequality (\ref{estnonhoms2}) is proved in
\cite{MR2101982}, proposition 2.8. $\Box$

\subsection{Technical lemmas}
In this subsection, we recall some useful lemmas which allow to
share derivatives of various expressions  in $L^p_xL^q_T$ norms.
One can
find proofs of lemmas \ref{leibniz}-\ref{lemg} in \cite{MR2101982,MR1086966}.\\
Here $f$ and $g$ denote two elements of $\mathcal{S}(\R)$.

\begin{lemma}\label{leibniz} If $\alpha>0$ and $1<p,q<\infty$, then
\[\|D^\alpha_x(fg)\|_{L^{p}_xL^{q}_T}\lesssim\|f\|_{L^{p_1}_xL^{q_1}_T}\|D^\alpha_xg\|_{L^{p_2}_xL^{q_2}_T}+
\|g\|_{L^{\tilde{p}_1}_xL^{\tilde{q}_1}_T}\|D^\alpha_xf\|_{L^{\tilde{p}_2}_xL^{\tilde{q}_2}_T}\]
where $1<p_1,p_2,q_2,\tilde{p}_1,\tilde{p}_2,\tilde{q}_2<\infty$,
$1<q_1,\tilde{q}_1\leq\infty$,
$1/p_1+1/p_2=1/\tilde{p}_1+1/\tilde{p}_2=1/p$ and
$1/q_1+1/q_2=1/\tilde{q}_1+1/\tilde{q}_2=1/q$.\\
Moreover the cases $(p_1,q_1)=(\infty,\infty)$ and
$(\tilde{p}_1,\tilde{q}_1)=(\infty,\infty)$ are allowed.
\end{lemma}

\begin{lemma}\label{comp} If $0<\alpha<1$ and $1<p,q<\infty$ then
\[\|D^\alpha_xF(f)\|_{L^p_xL^q_T}\lesssim
\|F'(f)\|_{L^{p_1}_xL^{q_1}_T}\|D^\alpha_xf\|_{L^{p_2}_xL^{q_2}_T}\]
where $1<p_1,p_2,q_2<\infty$, $1<q_1\leq\infty$, $1/p_1+1/p_2=1/p$
and $1/q_1+1/q_2=1/q$.
\end{lemma}

\begin{lemma}\label{commu} If $0<\alpha<1$, $0\leq\beta<1-\alpha$ and
$1<p,q<\infty$, then
\[\|D^\beta_x([D^\alpha_x,f]g)\|_{L^p_xL^q_T}\lesssim
\|g\|_{L^{p_1}_xL^{q_1}_T}\|D^{\alpha+\beta}_xf\|_{L^{p_2}_xL^{q_2}_T}\]
where $1<p_1,q_1,p_2,q_2<\infty$, $1/p_1+1/p_2=1/p$ and
$1/q_1+1/q_2=1/q$.\\ Moreover, if $\beta>0$ then $q_1=\infty$ is
allowed.
\end{lemma}

\begin{lemma}\label{P+P-} If $\alpha>0$, $\beta\geq 0$ and $1<p,q<\infty$ then
\[\|D^\alpha_xP_+(fP_-D^\beta_xg)\|_{L^p_xL^q_T}\lesssim
\|D^{\gamma_1}_xf\|_{L^{p_1}_xL^{q_1}_T}\|D^{\gamma_2}_xg\|_{L^{p_2}_xL^{q_2}_T}\]
where $1<p_1,q_1,p_2,q_2<\infty$, $1/p_1+1/p_2=1/p$,
$1/q_1+1/q_2=1/q$ and $\gamma_1\geq\alpha$,
$\gamma_1+\gamma_2=\alpha+\beta$.
\end{lemma}

As in \cite{MR2101982}, we introduce the bilinear operator $G$
defined by
$$G(f,g)=\mathcal{F}^{-1}\Big(\frac{1}{2}\int_\R\frac{\xi_1(\xi-\xi_1)}{i\xi}[\sgn(\xi_1)
+\sgn(\xi-\xi_1)]\hat{f}(\xi_1)\hat{g}(\xi-\xi_1)d\xi_1\Big).$$ We
easily verify that
\begin{equation}\label{defg1}G(f,f)=\partial^{-1}_x(f_x\mathcal{H}f_x)=\partial_x^{-1}(-i(P_+f_x)^2+i(P_-f_x)^2)\end{equation}
and
\begin{equation}\label{defg2}G(f,g)=\partial_x^{-1}(-iP_+f_xP_+g_x+iP_-f_xP_-g_x).\end{equation}

\begin{lemma}\label{lemg} If $0\leq\alpha\leq 1$ and $1<p,q<\infty$ then
$$\|D^\alpha_xG(f,g)\|_{L^p_xL^q_T}\lesssim
\|D^{\gamma_1}_xf\|_{L^{p_1}_xL^{q_1}_T}\|D^{\gamma_2}_xg\|_{L^{p_2}_xL^{q_2}_T}$$
where $0\leq\gamma_1,\gamma_2\leq 1$,
$\gamma_1+\gamma_2=\alpha+1$, $1<p_1,q_1,p_2,q_2<\infty$,
$1/p_1+1/p_2=1/p$ and $1/q_1+1/q_2=1/q$.
\end{lemma}

We will also need the following lemma in order to treat low
frequencies in the integral term.

\begin{lemma}\label{lowfreq} If $\alpha\geq 0$ and $1\leq p,q\leq\infty$ then
$$\|P_0(fD^\alpha_xg)\|_{L^p_xL^q_T}\lesssim
\|D^{\gamma_1}_xf\|_{L^{p_1}_xL^{q_1}_T}\|D^{\gamma_2}_xg\|_{L^{p_2}_xL^{q_2}_T}+\|P_0f\|_{L^{\tilde{p}_1}_x
L^{\tilde{q}_1}_T}\|D^\alpha_xP_0g\|_{L^{\tilde{p}_2}_xL^{\tilde{q}_2}_T}$$
where $\gamma_1,\gamma_2\geq 0$, $\alpha=\gamma_1+\gamma_2$,
$1<p_i,q_i,\tilde{p}_i,\tilde{q}_i<\infty$,
$1/p_1+1/p_2=1/\tilde{p}_1+1/\tilde{p}_2=1/p$ and
$1/q_1+1/q_2=1/\tilde{q}_1+1/\tilde{q}_2=1/q$.
\end{lemma}
Proof : We split the product $fD^\alpha_xg$ as follows :
\begin{equation}\label{decbf}fD^\alpha_xg=P_+fP_+D^\alpha_xg+ P_+fP_-D^\alpha_xg+
P_-fP_+D^\alpha_xg+P_-fP_-D^\alpha_xg.\end{equation} It is
sufficient to consider the contribution of the first two terms.
For the first one, we remark that
$$P_0[P_+fP_+(D^\alpha_xg)]=P_0[P_0(P_+f)P_0(P_+D^\alpha_xg)]$$
and thus using the continuity of $P_0$ on $L^p_xL^q_T$,
\begin{eqnarray*}\|P_0[P_+fP_+(D^\alpha_xg)]\|_{L^p_xL^q_T}
&\lesssim & \|P_0(P_+f)P_0(P_+D^\alpha_xg)\|_{L^p_xL^q_T}\\
&\lesssim &
\|P_0f\|_{L^{\tilde{p}_1}_xL^{\tilde{q}_1}_T}\|D^\alpha_xP_0g\|_{L^{\tilde{p}_2}_xL^{\tilde{q}_2}_T}.
\end{eqnarray*}
For the second term in (\ref{decbf}) we have typically
contributions of the form $P_0[P_0(P_+f)P_0(P_-D^\alpha_xg)]$
which are treated as above, and
$P_0[\Tilde{P}_+f\tilde{P}_-D^\alpha_xg]$. Using decomposition
(\ref{DLP}), one can write
\begin{eqnarray*}P_0(\tilde{P}_+f\tilde{P}_-D^\alpha_xg) &=&
P_0\Big(\sum_{j\in\Z}Q_j(\tilde{P}_+f)P_j(\tilde{P}_-D^\alpha_xg)+\sum_{j\in\Z}P_j(\tilde{P}_+f)Q_j(\tilde{P}_-D^\alpha_xg)\Big)\\
&& +P_0\big(\sum_{|p|\leq
2}\sum_{j\in\Z}Q_j(\tilde{P}_+f)Q_{k-j}(\tilde{P}_-D^\alpha_xg)\Big).
\end{eqnarray*}
By a careful analysis of the various localisations, we get
$$P_0(\Tilde{P}_+f\tilde{P}_-D^\alpha_xg) =
P_0\Big[\sum_{|p|\lesssim
1}\sum_{j\in\Z}Q_j(\tilde{P}_+f)Q_{j+p}(\Tilde{P}_-D^\alpha_xg)\Big].$$
Here we define the operators $Q_j^\lambda=2^{-\lambda
j}D^\lambda_xQ_j$. It follows that
\begin{eqnarray*}P_0(\Tilde{P}_+f\tilde{P}_-D^\alpha_xg)
=P_0\Big[\sum_{|p|\lesssim
1}\sum_{j\in\Z}Q_j^{-\gamma_1}(\tilde{P}_+D^{\gamma_1}_xf)Q^{\gamma_1}_{j+p}(\tilde{P}_-D^{\gamma_2}_xg)\Big].\end{eqnarray*}
Thus using Cauchy-Schwarz and H\"older inequalities, and
Littlewood-Paley theorem,
\begin{eqnarray*}\|P_0(\Tilde{P}_+f\tilde{P}_-D^\alpha_xg)\|_{L^p_xL^q_T}
&\lesssim & \Big\|\sum_{|p|\lesssim
1}\Big[\Big(\sum_{j\in\Z}|Q^{-\gamma_1}_j\tilde{P}_+D^{\gamma_1}_xf|^2\Big)^{1/2}
\Big(\sum_{j\in\Z}|Q^{\gamma_1}_j\tilde{P}_-D^{\gamma_2}_xg|^2\Big)^{1/2}\Big]\Big\|_{L^p_xL^q_T}\\
&\lesssim & \Big\|
\Big(\sum_{j\in\Z}|Q^{-\gamma_1}_j\tilde{P}_+D^{\gamma_1}_xf|^2\Big)^{\frac{1}{2}}\Big\|_{L^{p_1}_xL^{q_1}_T}
\Big\|\Big(\sum_{j\in\Z}|Q^{\gamma_1}_j\tilde{P}_-D^{\gamma_2}_xg|^2\Big)^{\frac{1}{2}}\Big\|_{L^{p_2}_xL^{q_2}_T}\\
&\lesssim &
\|D^{\gamma_1}_xf\|_{L^{p_1}_xL^{q_1}_T}\|D^{\gamma_2}_xg\|_{L^{p_2}_xL^{q_2}_T}.
\Box
\end{eqnarray*}

\section{Nonlinear estimates}
\subsection{Gauge transformation}

By a rescaling argument, it is sufficient to solve
\begin{equation}\label{resceq}u_t+\mathcal{H}u_{xx}=2u^ku_x\end{equation}
(equation with minus sign in front of the nonlinearity could be
treated in the same way). If $u\in\mathcal{C}([0,T];H^\infty(\R))$
is a smooth solution, we define the gauge
transformation\footnotemark[1]
\begin{equation}\label{gauge}w=P_+(e^{-iF}u),\quad
F=F(u)=\int_{-\infty}^xu^k(y,t)dy.\end{equation}
\footnotetext[1]{we can also set $F=\frac 12 \int_{-\infty}^xu^k$
in the non-rescaled case $u_t+\H u_{xx}=u^ku_x$.}

The rest of this subsection is devoted to the proof of the
following estimate.

\begin{proposition} Let be $k\geq 12$ and $s_k<s<1/2$. Let $u\in\mathcal{C}([0,T];H^\infty(\R))$ be a solution
of the Cauchy problem associated to (\ref{resceq}) with initial
data $u_0\in H^\infty(\R)$. Then there exist $\nu=\nu(s,k)>0$ and
a positive nondecreasing polynomial function $p_{k}$ such that
\begin{eqnarray}\nonumber\|u\|_{X^s_T}&\lesssim &\|u_0\|_{H^s}+T^\nu
p_{k}(\|u\|_{X^s_T})\|u\|_{X^s_T}\\
\label{estnonlin1} &&+(\|u_0\|_{H^s}^k+T^\nu
p_k(\|u\|_{X^s_T})\|u\|_{X^{s}_T})\|D^{s+1/2}_xw\|_{L^\infty_xL^2_T}.\end{eqnarray}
\end{proposition}
Proof : We start by splitting $u$ according to (\ref{smoothdec}).
Then, using that $|P_+u|=|P_-u|$ (since $u$ is real), we deduce
\begin{equation}\label{decu}\|u\|_{X^s_T}\lesssim
\|P_0u\|_{X^s_T}+\|\tilde{P}_+u\|_{X^s_T}.\end{equation} For the
low frequencies, we use the Duhamel formulation of (\ref{gBO}),
lemma \ref{lemestnonhom} and (\ref{estsg}) to get
\begin{eqnarray*}\|P_0u\|_{X^s_T} &\lesssim
& \|P_0u_0\|_{H^s}+\|P_0D^{1/2}_xu^{k+1}\|_{L^1_xL^2_T}\\
&\lesssim & \|u_0\|_{H^s}+\|u^{k+1}\|_{L^1_xL^2_T}\\ &\lesssim &
\|u_0\|_{H^s}+T^\nu\|u\|_{L^{k+1}_xL^{2k(k+1)}_T}^{k+1}\\
&\lesssim & \|u_0\|_{H^s}+T^\nu\|u\|_{X^s_T}.
\end{eqnarray*}

Now we consider the second term in the right-hand side of
(\ref{decu}). As mentioned in \cite{MR2101982}, $\tilde{P}_+u$
satisfies the dispersive equation
\[\partial_t(\tilde{P}_+u)+\H\partial^2_x(\tilde{P}_+u)=\tilde{P}_+(e^{iF}u^kw_x)-\tilde{P}_+(e^{iF}u^k\partial_x
P_-(e^{-iF}u)) +i\tilde{P}_+(u^{2k+1}).\] Thus, according to lemma
\ref{lemestnonhom} \begin{eqnarray*}\|\Tilde{P}u\|_{X^s_T}
&\lesssim &
\|V(t)u_0\|_{X^s_T}+\Big\|\int_0^tV(t-\tau)\tilde{P}_+(e^{iF}u^kw_x)(\tau)d\tau\Big\|_{X^s_T}\\
&\quad&+\|D^s_xu^{2k+1}\|_{L^{6/5}_{xT}}+
\|\tilde{P}_+(e^{iF}u^k\partial_x
P_-(e^{-iF}u))\|_{L^{(\frac{5}{6}+\frac{s}{3})^{-1}}_xL^{(\frac{5}{6}-\frac{2s}{3})^{-1}}_T}\\
& \lesssim & \|u_0\|_{H^s}+A+B+C.
\end{eqnarray*}
Obviously, \[B \lesssim
\|u^{2k}\|_{L^{3/2}_{xT}}\|D^s_xu\|_{L^6_{xT}} \lesssim
\|u\|_{L^{3k}_{xT}}^{2k}\|D^s_xu\|_{L^6_{xT}}\lesssim
T^\nu\|u\|_{X^s_T}^{2k+1}.\] Term $C$ has a structure
$P_+(fP_-g_x)$ thus by lemma \ref{P+P-}
\begin{eqnarray*}C &\lesssim &
\|D^{1/2}_x(e^{iF}u^k)\|_{L^{6/5}_xL^3_T}\|D^{1/2}_x(e^{-iF}u)\|_{L^{3/s}_xL^{(\frac{1}{2}-\frac{2s}{3})^{-1}}_T}\\
&\lesssim & C_1C_2.\end{eqnarray*} Using lemmas
\ref{leibniz}-\ref{comp}, we infer \begin{eqnarray}\nonumber C_1
&\lesssim &
\|D^{1/2}_xu^k\|_{L^{6/5}_xL^3_T}+\|D^{1/2}_xe^{iF}\|_{L^{(\frac{1}{2}-s)^{-1}}_xL^{2/s}_T}\|u^k\|_{L^{(\frac{1}{3}+s)^{-1}}_x
L^{(\frac{1}{3}-\frac{s}{2})^{-1}}_T}\\ \nonumber &\lesssim &
\|u\|_{L^{(k-1)(\frac{5}{6}-\frac{s}{3})^{-1}}_xL^{(k-1)(\frac{2s}{3}-\frac{1}{6})^{-1}}_T}^{k-1}\|D^{1/2}_xu\|_{L^{3/s}_xL^{(\frac{1}{2}-\frac{2s}{3}
)^{-1}}_T}\\ \nonumber &\quad &
+\|D^{-1/2}_x(u^ke^{iF})\|_{L^{(\frac{1}{2}-s)^{-1}}_xL^{2/s}_T}
\|u\|_{L^{k(\frac{1}{3}+s)^{-1}}_x
L^{k(\frac{1}{3}-\frac{s}{2})^{-1}}_T}^k\\ \nonumber &\lesssim &
T^\nu\|u\|_{X^s_T}^k+T^\nu\|u^k\|_{L^{1/(1-s)}_xL^{2/s}_T}\|u\|_{X^s_T}^k\\
\label{c1} &\lesssim &
T^\nu\|u\|_{X^s_T}^k+T^\nu\|u\|_{X^s_T}^{2k}
\end{eqnarray}
and in the same way
\begin{eqnarray}\nonumber C_2 &\lesssim &\|D^{1/2}_xu\|_{L^{3/s}_xL^{(\frac{1}{2}-\frac{2s}{3})^{-1}}_T}+
\|D^{1/2}_xe^{-iF}\|_{L^{(\frac{4s}{3}-\frac{1}{2})^{-1}}_x
L^{(\frac{1}{2}-\frac{2s}{3})^{-1}}_T}\|u\|_{L^{(\frac{1}{2}-s)^{-1}}_xL^\infty_T}\\
\nonumber &\lesssim
&\|D^{1/2}_xu\|_{L^{3/s}_xL^{(\frac{1}{2}-\frac{2s}{3})^{-1}}_T}+\|u^k\|_{L^{3/4s}_x
L^{(\frac{1}{2}-\frac{2s}{3})^{-1}}_T}\|u\|_{L^{(\frac{1}{2}-s)^{-1}}_xL^\infty_T}\\
\label{c2} &\lesssim & \|u\|_{X^s_T}+T^\nu\|u\|_{X^s_T}^{k+1}.
\end{eqnarray}
Combining (\ref{c1}) and (\ref{c2}), $C$ is bounded by $$C\lesssim
T^\nu(\|u\|_{X^s_T}^{k+1}+\|u\|_{X^s_T}^{2k+1}+\|u\|_{X^s_T}^{3k+1})\lesssim
T^\nu p_k(\|u\|_{X^s_T})\|u\|_{X^s_T}.$$ In order to study the
contribution of $A$, we decompose $e^{iF}u^kw_x$ as
\[e^{iF}u^kw_x=D^{1/2}_x(e^{iF}u^k\H D^{1/2}_xw)-[D^{1/2}_x,e^{iF}u^k]\H D^{1/2}_xw.\]
Therefore, according to lemma \ref{lemestnonhom}, and using the
fact that $\tilde{P}_+$ is continuous on $L^1_xL^2_T$,
\begin{eqnarray*}A &\lesssim &
\|D^s_x(e^{iF}u^k\H
D^{1/2}_xw)\|_{L^1_xL^2_T}+\|[D^{1/2}_x,e^{iF}u^k]\H
D^{1/2}_xw\|_{L^{(\frac{5}{6}+\frac{s}{3})^{-1}}_x
L^{(\frac{5}{6}-\frac{2s}{3})^{-1}}_T}\\ &\lesssim & A_1+A_2.
\end{eqnarray*}
Note that $A_1$ cannot be treated by lemma \ref{leibniz}, so we
use lemma A.13 in \cite{MR1086966}. This leads to
\begin{eqnarray*}A_1 &\lesssim &
\|D^s_x(e^{iF}u^k)\|_{L^{(1-\frac{s}{3})^{-1}}_xL^{3/2s}_T}
\|D^{1/2}_x(e^{-iF}u)\|_{L^{3/s}_xL^{(\frac{1}{2}-\frac{2s}{3})^{-1}}_T}+\|u^k\|_{L^1_xL^\infty_T}\|D^{s+1/2}_xw\|_{L^\infty_xL^2_T}\\
&\lesssim & A_{11}C_2+A_{12}^k\|D^{s+1/2}_xw\|_{L^\infty_xL^2_T}.
\end{eqnarray*}
By lemma \ref{leibniz} we bound the contribution of $A_{11}$ by
\begin{eqnarray*}
A_{11} &\lesssim &
\|D^s_xu^k\|_{L^{(1-\frac{s}{3})^{-1}}_xL^{3/2s}_T}+\|D^s_xe^{iF}\|_{L^{3/2s}_xL^{6/s}_T}
\|u^k\|_{L^{1/(1-s)}_xL^{2/s}_T}\\ &\lesssim &
\|u\|_{L^{(k-1)(\frac{5}{6}-\frac{s}{3})^{-1}}_xL^{(k-1)(\frac{2s}{3}-\frac{1}{6})^{-1}}_T}^{k-1}\|D^s_xu\|_{L^6_{xT}}+
\|u\|_{L^{k(1-\frac{s}{3})^{-1}}_xL^{6k/s}_T}^k
\|u\|_{L^{k/(1-s)}_xL^{2k/s}_T}^k \\ &\lesssim &
T^\nu\|u\|_{X^s_T}^k+T^\nu\|u\|_{X^s_T}^{2k}.
\end{eqnarray*}
To treat $A_{12}=\|u\|_{L^k_xL^\infty_T}$ we use the Duhamel
formulation of (\ref{gBO}) and lemma \ref{estsnonhom},
\begin{eqnarray*}A_{12} &\lesssim &
\|V(t)u_0\|_{L^k_xL^\infty_T}+\Big\|\int_0^tV(t-\tau)\partial_xu^{k+1}(\tau)d\tau\Big\|_{L^k_xL^\infty_T}\\
&\lesssim &
\|u_0\|_{H^s}+\|D^{s_k+1/2+3\eps}_xu^{k+1}\|_{L^{(1-\eps)^{-1}}_xL^{(\frac
12+2\eps)^{-1}}_T}
\end{eqnarray*}
and setting $\eps'=\frac 13(s-s_k)-\eps>0$ it follows that
\begin{eqnarray*}
\lefteqn{\|D^{s_k+1/2+3\eps}_xu^{k+1}\|_{L^{(1-\eps)^{-1}}_xL^{(\frac
12+2\eps)^{-1}}_T}}\\ &\lesssim &
\|D^{s+1/2-3\eps'}_xu\|_{L^{1/\eps'}_xL^{(\frac
12-2\eps')^{-1}}_T} \|u^k\|_{L^{(1-\frac
13(s-s_k))^{-1}}_xL^{\frac 32(s-s_k)^{-1}}_T}\\ &\lesssim &
T^\nu\|u\|_{X^s_T}\|u\|_{L^{k(1-\frac
13(s-s_k))^{-1}}_xL^\infty_T}^k
\\ &\lesssim &
T^\nu\|u\|_{X^s_T}^{k+1}.
\end{eqnarray*}
Finally, according to lemma \ref{commu} we write
\begin{eqnarray*}A_2 &\lesssim &
\|D^{1/2}_x(e^{iF}u^k)\|_{L^{6/5}_xL^3_T}
\|D^{1/2}_x(e^{-iF}u)\|_{L^{3/s}_xL^{(\frac{1}{2}-\frac{2s}{3})^{-1}}_T}\\
&\lesssim & C_1C_2\\ &\lesssim & T^\nu
p_k(\|u\|_{X^s_T})\|u\|_{X^s_T},
\end{eqnarray*}
witch complete the proof of (\ref{estnonlin1}). $\Box$

\subsection{Estimate of $\|D^{s+1/2}_xw\|_{L^\infty_xL^2_T}$}

Now our aim is to estimate the term
$\|D^{s+1/2}_xw\|_{L^\infty_xL^2_T}$ which appears in
(\ref{estnonlin1}). More precisely we will prove the following
proposition.

\begin{proposition} Let $k\geq 12$ and $s_k<s<1/2$. For all solution $u\in\mathcal{C}([0,T];
H^\infty(\R))$ of (\ref{resceq}) with initial data
$u_0\in{H}^\infty(\R)$, we have the following bound,
\begin{equation}\label{estnonlin2}\|D^{s+1/2}_xw\|_{L^\infty_xL^2_T} \lesssim p_{k}(\|u_0\|_{H^s})\|u_0\|_{H^s}+T^\nu
p_{k}(\|u\|_{X^s_T})\|u\|_{X^s_T}\end{equation} where $p_{k}$ is a
positive nondecreasing polynomial function.
\end{proposition}
Proof : Following \cite{MR2101982}, we see that $w$ satisfies the
equation \begin{eqnarray}\nonumber w_t+\H w_{xx} &=&
P_+[2e^{-iF}(-ku^kP_-u_x-iP_-u_{xx})]\\
\label{eqw} &&-ik(k-1)P_+\Big(e^{-iF}u\int_{-\infty}^xu^{k-2}u_x\H
u_x\Big).\end{eqnarray} Thus using the Duhamel formulation of
(\ref{eqw}) and lemma \ref{lemestnonhom} we infer
\begin{eqnarray}\nonumber \|D^{s+1/2}_xw\|_{L^\infty_xL^2_T} &\lesssim &
\|D^{s+1/2}_xV(t)w(0)\|_{L^\infty_xL^2_T}\\ &&\nonumber +
\|P_+[2e^{-iF}(-ku^kP_-u_x-iP_-u_{xx})]\|_{L^{(\frac{5}{6}+\frac{s}{3})^{-1}}_xL^{(\frac{5}{6}-\frac{2s}{3})^{-1}}_T}\\
\label{termint}&\quad &
+\Big\|D^{s+1/2}_x\int_0^tV(t-\tau)P_+\Big(e^{-iF}u\int_{-\infty}^xu^{k-2}u_x\H
u_x\Big)d\tau\Big\|_{L^\infty_xL^2_T}.
\end{eqnarray}
The first term of right-hand side can be bounded by
\begin{eqnarray*}\|D^{s+1/2}_xV(t)w(0)\|_{L^\infty_xL^2_T} &\lesssim &
\|e^{-iF(u_0)}u_0\|_{H^s}\\ &\lesssim &
\|u_0\|_{L^2}+\|e^{-iF(u_0)}\|_{L^\infty}\|D^s_xu_0\|_{L^2}+\|D^s_xe^{-iF(u_0)}\|_{L^{1/s}}\|u_0\|_{L^{(\frac{1}{2}-s)^{-1}}}\\
&\lesssim &
\|u_0\|_{H^s}+\|D^{s-1}(e^{-iF(u_0)}u_0^k)\|_{L^{1/s}}\|u_0\|_{H^s}\\
&\lesssim & \|u_0\|_{H^s}(1+\|u_0^k\|_{L^1})\\
&\lesssim & \|u_0\|_{H^s}(1+\|u_0\|_{H^s}^k).\end{eqnarray*} On
the other hand, according to lemma \ref{P+P-}, we see that
\begin{eqnarray*} \lefteqn{\|P_+[2e^{-iF}(-ku^kP_-u_x-iP_-u_{xx})]\|_{L^{(\frac{5}{6}+\frac{s}{3})^{-1}}_xL^{(\frac{5}{6}-\frac{2s}{3})^{-1}}_T}}
\\ &\lesssim &
\|D^{1/2}_x(e^{-iF}u^k)\|_{L^{6/5}_xL^3_T}\|D^{1/2}_xu\|_{L^{3/s}_xL^{(\frac{1}{2}-\frac{2s}{3})^{-1}}_T}\\
&\lesssim & C_1\|u\|_{X^s_T}\\ &\lesssim &
T^\nu\|u\|_{X^s_T}^{k+1}+T^\nu\|u\|_{X^s_T}^{2k+1}
\end{eqnarray*}
Thus it remains to estimate the integral term in (\ref{termint}),
that is, the last one. For this purpose, we split it as
\begin{eqnarray*}\int_{-\infty}^xu^{k-2}u_x\H u_x &=& {P}_0
\int_{-\infty}^xu^{k-2}u_x\H u_x+
\tilde{P}_+\int_{-\infty}^xu^{k-2}u_x\H
u_x+\tilde{P}_-\int_{-\infty}^xu^{k-2}u_x\H u_x\\
&=& I+II+III.\end{eqnarray*} By symmetry, it will be enough to
consider the contributions of $I$ and $II$.

\textbf{Contribution of $I$}\\
Using a commutator operator, we decompose \[D^s_x(e^{-iF}u
I)=D^s_x(e^{-iF}u)I+[D^s_x,I]e^{-iF}u.\] Therefore thanks to lemma
\ref{lemestnonhom} we obtain
\begin{eqnarray*}\lefteqn{\Big\|D^{s+1/2}_x \int_{0}^tV(t-\tau)P_+\Big(e^{-iF}u\int_{-\infty}^x{P}_0(u^{k-2}u_x\H
u_x)\Big)d\tau\Big\|_{L^\infty_xL^2_T}} \\ & \lesssim &
\Big\|D^s_x(e^{-iF}u)\int_{-\infty}^x{P}_0(u^{k-2}u_x\H
u_x)\Big\|_{L^1_TL^2_x}  \\ && \quad  +
\Big\|D^{1/4}_x\Big[D^s_x,\int_{-\infty}^x{P}_0(u^{k-2}u_x\H
u_x)\Big]e^{-iF}u\Big\|_{L^{4/3}_xL^1_T}\\ &  \lesssim &
D+E.\end{eqnarray*} The contribution of $D$ is treated as follows.
\begin{eqnarray*}D &\lesssim &
\|D^s_x(e^{-iF}u)\|_{L^\infty_TL^2_x}\Big\|\int_{-\infty}^x{P}_0(u^{k-2}u_x\H
u_x)\Big\|_{L^1_TL^\infty_x}\\ &\lesssim &
(\|D^s_xu\|_{L^\infty_TL^2_x}+\|D^s_xe^{-iF}\|_{L^\infty_TL^{1/s}_x}\|u\|_{L^\infty_TL^{(\frac{1}{2}-s)^{-1}}_x})
\|{P}_0(u^{k-2}\partial_xG(u,u))\|_{L^1_{xT}}\\
&\lesssim &
T^\nu\|u\|_{L^\infty_TH^s_x}(1+\|u\|_{L^\infty_TH^{s_k}_x}^k)
\|{P}_0(u^{k-2}\partial_xG(u,u))\|_{L^1_xL^{1/(1-\eps)}_{T}}\\
&\lesssim &
T^\nu\|u\|_{X^s_T}(1+\|u\|_{X^s_T}^k)\|{P}_0(u^{k-2}\partial_xG(u,u))\|_{L^1_xL^{1/(1-\eps)}_{T}}.
\end{eqnarray*}
The low frequencies term is estimated with lemma \ref{lowfreq}. We
get
\begin{eqnarray*}
\lefteqn{\|{P}_0(u^{k-2}\partial_xG(u,u))\|_{L^1_xL^{1/(1-\eps)}_{T}}}
\\ &\lesssim &
\|D^{1-2s+6\eps}_xu^{k-2}\|_{L^{1/(1-2\eps)}_xL^{1/3\eps}_T}
\|D^{2s-6\eps}_xG(u,u)\|_{L^{1/2\eps}_xL^{1/(1-4\eps)}_T} \\
&\quad &
+\|{P}_0u^{k-2}\|_{L^{(1-\frac{2s}{3})^{-1}}_xL^{(\frac{4s}{3}-\eps)^{-1}}_T}\|{P}_0\partial_x
G(u,u)\|_{L^{3/2s}_xL^{(1-\frac{4s}{3})^{-1}}_T}\\ &\lesssim &
\|u^{k-3}\|_{L^{(3s-\frac{1}{2}-2\eps)^{-1}}_xL^{\infty}_T}
\|D^{1-2s+6\eps}_xu\|_{L^{(\frac{3}{2}-3s)^{-1}}_xL^{1/3\eps}_T}
\|D^{s+1/2-3\eps}_xu\|_{L^{1/\eps}_xL^{(\frac{1}{2}-2\eps)^{-1}}_T}^2\\
&\quad & +\|u\|_{L^{k-2}_xL^\infty_T}^{k-2}
\|D^{1/2}_xu\|_{L^{3/s}_xL^{(\frac{1}{2}-\frac{2s}{3})^{-1}}_T}^2\\
&\lesssim & \|u\|_{X^s_T}^{k}.
\end{eqnarray*}
Note that in order to bound the norm
$N_9=\|D^{1-2s+6\eps}_xu\|_{L^{(\frac{3}{2}-3s)^{-1}}_xL^{1/3\eps}_T}$,
we have to impose $k\geq 12$. Indeed, for $\eps>0$ small enough,
the triplet $(1-3s+6\eps,(\frac 32-3s)^{-1},1/3\eps)$ is
1-admissible if and only if $$\Big(\frac 32-3s\Big)^{-1}\geq
4\quad \textrm{ and }\quad 2\Big(\frac 32-3s\Big)+3\eps\leq \frac
12$$ if and only if $s>5/12=1/2-1/12$.

To bound $E$ by lemma \ref{commu},
\begin{eqnarray*}E &\lesssim & T^\nu \Big\|D^{1/4}_x\Big[D^s_x,\int_{-\infty}^x{P}_0(u^{k-2}u_x\H
u_x)\Big]e^{-iF}u\Big\|_{L^{4/3}_xL^{1/(1-\eps)}_T}\\
&\lesssim & T^\nu
\|D^{s-3/4}_x{P}_0(u^{k-2}\partial_xG(u,u))\|_{L^{(s+\frac{1}{4})^{-1}}_xL^{1/(1-\eps)}_T}
\|u\|_{L^{(\frac{1}{2}-s)^{-1}}_xL^\infty_T}\\ &\lesssim & T^\nu
\|{P}_0(u^{k-2}\partial_xG(u,u))\|_{L^1_xL^{1/(1-\eps)}_T}\|u\|_{X^s_T}\\
&\lesssim & T^\nu\|u\|_{X^s_T}^{k+1}.
\end{eqnarray*}

\textbf{Contribution of $II$}\\
We split the term $II$ into $$II = II_1+II_2+II_3$$ with
\begin{eqnarray*} II_1 &=&
\int_{-\infty}^x\tilde{P}_+(u^{k-2}P_-(u_x\H u_x)),\\ II_2 &=&
\int_{-\infty}^x\tilde{P}_+(P_+(u^{k-2})P_+(u_x\H u_x)),\\
II_3&= & \int_{-\infty}^x\tilde{P}_+(P_-(u^{k-2})P_+(u_x\H u_x)).
\end{eqnarray*}

\textbf{Contribution of $II_1$}\\
The treatment of $II_1$ is similar to the one of $I$. We write
\[D^s_x(e^{-iF}u II_1)=D^s_x(e^{-iF}u)II_1+[D^s_x,II_1]e^{-iF}u\]
and thus
\begin{eqnarray*}\lefteqn{\Big\|D^{s+1/2}_x\int_{0}^tV(t-\tau)P_+\Big(e^{-iF}u\int_{-\infty}^x\tilde{P}_+(u^{k-2}P_-(u_x\H
u_x))\Big)d\tau\Big\|_{L^\infty_xL^2_T}} \\ &\lesssim &
\Big\|D^s_x(e^{-iF}u)\int_{-\infty}^x\tilde{P}_+(u^{k-2}P_-(u_x\H
u_x))\Big\|_{L^1_TL^2_x} \\ && +
\Big\|D^{1/4}_x\Big[D^s_x,\int_{-\infty}^x\tilde{P}_+(u^{k-2}P_-(u_x\H
u_x))\Big]e^{-iF}u\Big\|_{L^{4/3}_xL^1_T}\\&\lesssim &
D'+E'.\end{eqnarray*} We first bound $D'$ as
\begin{eqnarray*}D' &\lesssim &
\|D^s_x(e^{-iF}u)\|_{L^\infty_TL^2_x}\Big\|\int_{-\infty}^x\tilde{P}_+(u^{k-2}P_-(u_x\H
u_x))\Big\|_{L^1_TL^\infty_x}\\ &\lesssim &
\|u\|_{X^s_T}(1+\|u\|_{X^s_T}^k)
\|\tilde{P}_+(u^{k-2}P_-\partial_xG(u,u))\|_{L^1_TL^{1/(1-\eps)}_{x}}\\
&\lesssim & T^\nu\|u\|_{X^s_T}(1+\|u\|_{X^s_T}^k)
\|\tilde{P}_+(u^{k-2}P_-\partial_xG(u,u))\|_{L^{1/(1-\eps)}_{xT}}
\end{eqnarray*}
and using lemma \ref{P+P-}, we get
\begin{eqnarray*}
\lefteqn{\|\tilde{P}_+(u^{k-2}P_-\partial_xG(u,u))\|_{L^{1/(1-\eps)}_{xT}}}
\\ &\lesssim &
\|D^{1-2s+6\eps}_xu^{k-2}\|_{L^{1/(1-3\eps)}_xL^{1/3\eps}_T}
\|D^{2s-6\eps}_xG(u,u)\|_{L^{1/2\eps}_xL^{1/(1-4\eps)}_T}\\
&\lesssim &
\|u^{k-3}\|_{L^{(3s-\frac{1}{2}-3\eps)^{-1}}_xL^{\infty}_T}
\|D^{1-2s+6\eps}_xu\|_{L^{(\frac{3}{2}-3s){-1}}_xL^{1/3\eps}_T}
\|D^{s+1/2-3\eps}_xu\|_{L^{1/\eps}_xL^{(\frac{1}{2}-2\eps)^{-1}}_T}^2
\\
&\lesssim & \|u\|_{X^s_T}^{k}.
\end{eqnarray*}
Next, $E'$ is estimated as follows
\begin{eqnarray*}E' &\lesssim & T^\nu \Big\|D^{1/4}_x\Big[D^s_x,\int_{-\infty}^x\tilde{P}_+(u^{k-2}P_-(u_x\H
u_x))\Big]e^{-iF}u\Big\|_{L^{4/3}_xL^{1/(1-\eps)}_T}\\
&\lesssim & T^\nu
\|D^{s-3/4}_x\tilde{P}_+(u^{k-2}P_-\partial_xG(u,u))\|_{L^{(s+\frac{1}{4})^{-1}}_xL^{1/(1-\eps)}_T}
\|u\|_{L^{(\frac{1}{2}-s)^{-1}}_xL^\infty_T}\\ &\lesssim & T^\nu
\|\tilde{P}_+(u^{k-2}P_-\partial_xG(u,u))\|_{L^{1/(1-\eps)}_{xT}}\|u\|_{X^s_T}\\
&\lesssim & T^\nu\|u\|_{X^s_T}^{k+1}.
\end{eqnarray*}

\textbf{Contribution of $II_2$}\\
A decomposition of $u^{k-2}$ into low and high frequencies, an
integration by parts, and formulas (\ref{defg1}-\ref{defg2}) give
\begin{equation}\begin{split} II_2 &= \int_{-\infty}^x\tilde{P}_+[P_+P_{\leq
-4}u^{k-2}P_+P_{\geq
-3}\partial_xG(u,u)]+\int_{-\infty}^x\tilde{P}_+[P_+P_{\geq
-3}u^{k-2}P_+\partial_xG(u,u)]\\ \label{II2} &=
\tilde{P}_+[P_+P_{\leq -4}u^{k-2}P_+P_{\geq
-3}G(u,u)]-i\tilde{P}_+G(P_{\leq -4}u^{k-2}, P_{\geq
-3}\partial_x^{-1}G(u,u))\\  & \quad +i \tilde{P}_+G(P_{\geq
-3}\partial_x^{-1}u^{k-2},G(u,u)).\end{split}
\end{equation}
We bound the first term by
\begin{eqnarray}\nonumber \lefteqn{\Big\| D^{s+1/2}_x\int_0^tV(t-\tau)P_+(e^{-iF}u\tilde{P}_+
[P_+P_{\leq -4}u^{k-2}P_+P_{\geq
-3}G(u,u)])(\tau)d\tau\Big\|_{X^s_T}}
\\ \nonumber &\lesssim &
\|u\tilde{P}_+[P_+P_{\leq -4}u^{k-2}P_+P_{\geq
-3}G(u,u)]\|_{L^{(\frac{5}{6}+\frac{s}{3})^{-1}}_xL^{(\frac{5}{6}-\frac{2s}{3})^{-1}}_T}\\
\nonumber &\lesssim &
\|u\|_{L^{(k-1)(\frac{5}{6}-\frac{s}{3})^{-1}}_xL^{(k-1)(\frac{2s}{3}-\frac{1}{6})^{-1}}_T}\\
\nonumber &&\times \|P_+P_{\leq -4}u^{k-2}P_+P_{\geq
-3}G(u,u)\|_{L^{(k-1)[k(\frac{5}{6}+\frac{s}{3})-\frac{5}{3}]^{-1}}_xL^{(k-1)[k(\frac{5}{6}-\frac{2s}{3})-\frac{2}{3}]^{-1}}_T}
\\
 &\lesssim & \nonumber
\|u\|_{L^{(k-1)(\frac{5}{6}-\frac{s}{3})^{-1}}_xL^{(k-1)(\frac{2s}{3}-\frac{1}{6})^{-1}}_T}\\
\nonumber && \times
\|u^{k-2}\|_{L^{\frac{k-1}{k-2}(\frac{5}{6}-\frac{s}{3})^{-1}}_xL^{\frac{k-1}{k-2}(\frac{2s}{3}-\frac{1}{6})^{-1}}_T}
\|G(u,u)\|_{L^{3/2s}_xL^{(1-\frac{4s}{3})^{-1}}_T}\\
\nonumber &\lesssim &
\|u\|_{L^{(k-1)(\frac{5}{6}-\frac{s}{3})^{-1}}_xL^{(k-1)(\frac{2s}{3}-\frac{1}{6})^{-1}}_T}^{k-1}
\|D^{1/2}_xu\|_{L^{3/s}_xL^{(\frac{1}{2}-\frac{2s}{3})^{-1}}_T}^2\\
\label{II21} &\lesssim &  T^\nu\|u\|_{X^s_T}^{k+1}.
\end{eqnarray}
The other terms in (\ref{II2}) are teated in the same way via
lemma \ref{lemg}.

\textbf{Contribution of $II_3$}\\
In order to share the derivative on $G(u,u)$ in $II_3$ with lemma
\ref{P+P-}, we first integrate by parts
\[II_3 =
\tilde{P}_+[P_-u^{k-2}P_+G(u,u)]-\tilde{P}_+\Big[\int_{-\infty}^x
P_-\partial_x u^{k-2}P_+G(u,u)\Big].\] Then we see that the first
term can be estimated exactly as (\ref{II21}). Finally for the
last term in the previous equality we repeat the proof for the
contribution of $II_1$. $\Box$

\section{Proof of Theorem \ref{mainresult}}
In this section we briefly recall the standard arguments which
yield well-posedness for (\ref{gBO}) ; we refer the reader to
\cite{MR2101982} for details. We choose $k\geq 12$ and
$s_k<s<1/2$.

We start by taking a sequence $(u_0^n)_n$ in $H^\infty(\R)$ such
that $u_0^n\rightarrow u_0$ in $H^s(\R)$ and $\|u_0^n\|_{H^s}\leq
\|u_0\|_{H^s}$. Now let $u_n\in H^\infty(\R)$ be the solutions of
(\ref{gBO}) with initial data $u_0^n$. Then bounds
(\ref{estnonlin1}) and (\ref{estnonlin2}) imply the \textit{a
priori} estimate
\begin{equation}\|u_n\|_{X^s_T}\lesssim
p_{k}(\|u_0^n\|_{H^s})\|u_0^n\|_{H^s}+T^\nu
p_{k}(\|u_n\|_{X^s_T})\|u_n\|_{X^s_T}.\end{equation} This allows
us to obtain the existence of a $T>0$ small enough and a solution
$u\in X^s_T$ of (\ref{gBO}).

Using the integral equation (\ref{duhamel}) and
(\ref{estnonlin1})-(\ref{estnonlin2}) it follows that for all
$0<t_1<t_2<T$,
\begin{eqnarray*}\|u(t_1)-u(t_2)\|_{H^s} &\lesssim &
\sup_{t\in[t_1,t_2]}\|u(t)-u(t_1)\| \\ &\lesssim &
\|u(t)-u(t_1)\|_{L^\infty([t_1,t_2];H^s)}+\Big\|\int_{t_1}^tV(t-\tau)\partial_x(u^{k+1})(\tau)d\tau\Big\|_{L^\infty([t_1,t_2];H^s)}
\\ &\lesssim & o(1).
\end{eqnarray*}
 This shows that $u\in\mathcal{C}([0,T];H^s(\R))$.

  We now turn to the proof of the uniqueness and the dependance of
 the solution upon the data. In this purpose we must establish the
 estimate
 $$\|u_1-u_2\|_{X^s_T}\lesssim
 p_{k}(\|u_{0,1}\|_{H^s}+\|u_{0,2}\|_{H^s})\|u_{0,1}-u_{0,2}\|_{H^s}+T^\nu
 p_{k}(\|u_1\|_{X^s_T}+\|u_2\|_{X^s_T})\|u_1-u_2\|_{X^s_T}$$ for
 $u_1$, $u_2$ two solutions of (\ref{gBO}) associated to initial
 data $u_{0,1}$ and $u_{0,2}$ respectively. We process exactly as
 in section 3 with the gauge transformation $$w=w_1-w_2,\quad
 w_j=P_+(e^{-iF_j}u_j),\quad F_j=\int_{-\infty}^xu^k_j(y,t)dy.$$
 The main new ingredient to use is the estimate
 $$|e^{i\int_{-\infty}^xf_1}-e^{i\int_{-\infty}^xf_2}|\lesssim
 \|f_1-f_2\|_{L^1}$$ for any real functions $f_1$, $f_2$ as explained in \cite{MR2101982}.

  \renewcommand{\theequation}{A-\arabic{equation}}
  % redefine the command that creates the equation no.
  \setcounter{equation}{0}  % reset counter
  \section*{Appendix}  % use *-form to suppress numbering

This subsection is devoted to the proof of theorem \ref{illth}. As
in \cite{MR2038121,MR1950814,MR1885293}, it is a consequence of
the following result.

\begin{lemma}\label{lem}Let $s<1/3$. Then there exists a sequence of
functions $\{h_N\}\subset H^s(\mathbb{R})$ such that for all
$T>0$,
\[\|h_N\|_{H^s}\lesssim 1\]
\begin{equation}\label{cond2}\lim_{N\rightarrow
+\infty}\sup_{[0,T]}\Big\|\int_0^tV(t-s)\partial_x((V(s)h_N)^4)ds\Big\|_{H^s}=+\infty\end{equation}

\end{lemma}

We show first that lemma \ref{lem} implies the result.\\
Suppose that theorem \ref{illth} fails. Since the flow-map
$\varphi\mapsto u(\varphi)$ is of class $\mathcal{C}^4$ at the
origin, we have the relation
\[F(u,\varphi):=u(\varphi)-V(t)\varphi+\int_0^tV(t-s)u^3(s)\partial_xu(s)ds=0\]
which together with the implicit function theorem yields
\[v(t,x):=\frac{\partial^3
u}{\partial\varphi^3}(t,x,0)[h_N,...,h_N]=3!\int_0^tV(t-s)\partial_x((V(s)h_N)^4)ds.\]
Hence \[\sup_{[0,T]}\|v(t)\|_{H^s}\lesssim\|h_N\|_{H^s}\lesssim
1,\] witch contradicts (\ref{cond2}).

\textbf{Proof of lemma \ref{lem}}\\
For each integer $N$, we define the function $h_N$ though its
Fourier transform by
\[\widehat{h_N}(\xi)=\alpha^{-1/2}N^{-s}\big(\chi_1(\xi)+\chi_2(\xi)\big)\]
where $\chi_1=\chi_{[N,N+\alpha]}$, $\chi_2(\xi)=\chi_1(-\xi)$ and
$\alpha=N^{-\theta}$, $\theta>0$ to be chosen later. Observe that
$h_N$ is a real valued function since $\widehat{h_N}$ is even.
Moreover, an obvious calculation yields $\|h_N\|_{H^s}\simeq 1$.\\
We now want to estimate $|\hat{v}|$ :
\begin{eqnarray*}\hat{v}(\xi_0,t) & \simeq & \xi_0 e^{ip(\xi_0)t}\int_0^te^{-ip(\xi_0)s}\mathcal{F}_x((V(s)h_N)^4)ds\\
 &\simeq & \alpha^{-2}N^{-4s}\xi_0
 e^{ip(\xi_0)t}\sum_{n=0}^4\int_0^te^{-ip(\xi_0)s}\big(e^{ip(\xi_0)s}\chi_1\big)^{\ast
 n}\ast\big(e^{ip(\xi_0)s}\chi_2\big)^{\ast(4-n)}ds\\ & := &
 \sum_{n=0}^4v_n
 \end{eqnarray*}
 where we defined $p(\xi)=\xi|\xi|$ and $f^{\ast n}=f\ast...\ast f$.
The function $v_4$ is supported in $[4N,4N+4\alpha]$ which is
disjoined with the supports of $v_n$, $n=0,1,2,3$. Consequently,
\begin{eqnarray*}\lefteqn{\hat{v}(\xi_0,t)\chi_{[4N,4N+4\alpha]}(\xi_0)}\\ &
\simeq & \alpha^{-2}N^{-4s}\xi_0
e^{ip(\xi_0)t}\int_0^t\int_{\mathbb{R}^3}
e^{-ip(\xi_0)s}e^{ip(\xi_0-\xi_1)s}e^{ip(\xi_1-\xi_2)s}
e^{ip(\xi_2-\xi_3)s}e^{ip(\xi_3)s}\\
 & &  \quad \times \chi_1(\xi_0-\xi_1)\chi_1(\xi_1-\xi_2)\chi_1(\xi_2-\xi_3)\chi_1(\xi_3)d\xi_1d\xi_2d\xi_3ds\\
& \simeq & \alpha^{-2}N^{-4s}\xi_0
e^{ip(\xi_0)t}\int_{\mathbb{R}^3}\frac{e^{itP(\xi_0,\xi_1,\xi_2,\xi_3)}-1}{P(\xi_0,\xi_1,\xi_2,\xi_3)}\\
& & \quad \times
\chi_1(\xi_0-\xi_1)\chi_1(\xi_1-\xi_2)\chi_1(\xi_2-\xi_3)\chi_1(\xi_3)d\xi_1d\xi_2d\xi_3
\end{eqnarray*}
with
$P(\xi_0,\xi_1,\xi_2,\xi_3)=-2\sum_{j=1}^3\xi_j(\xi_{j-1}-\xi_j)$.
For $(\xi_{j-1}-\xi_j)$ and $\xi_3$ in $[4N,4N+4\alpha]$, we have
\[P(\xi_0,\xi_1,\xi_2,\xi_3)\simeq N^2\quad\textrm{ and }\quad\Big|\frac{e^{itP(\xi_0,\xi_1,\xi_2,\xi_3)}-1}{P(\xi_0,\xi_1,\xi_2,\xi_3)}\Big|
=|t|+O(N^2) \gtrsim 1,\] witch yields
\[|\hat{v}(\xi_0,t)|\chi_{[4N,4N+4\alpha]}(\xi_0)\gtrsim
\alpha^{-2}N^{-4s}|\xi_0|\chi_1^{\ast4}(\xi_0).\] By
straightforward calculations,
\[\chi_1^{\ast4}(\xi)\simeq\int_{\mathbb{R}}\exp\big(-i\xi_1(4N+2\alpha-\xi)\big)\Big(\frac{\sin
i\alpha\xi_1/2}{\xi_1}\Big)^4d\xi_1\] and hence
$\chi_1^{\ast4}(4N+2\alpha)\simeq \alpha^3$. By a continuity
argument, $\chi_1^{\ast4}(\xi)\simeq \alpha^3$ for all
$\xi\in[4N,4N+4\alpha]$. This proves that
\[|\hat{v}(\xi_0,t)|\chi_{[4N,4N+4\alpha]}(\xi_0)\gtrsim \alpha
N^{-4s+1}\chi_{[4N,4N+4\alpha]}(\xi_0)\] and finally
\[\|v\|_{H^s}\gtrsim \alpha
N^{-4s+1}\Big(\int_{4N}^{4N+4\alpha}(1+|\xi|^2)^sd\xi\Big)^{1/2}\gtrsim
\alpha N^{-4s+1}N^s\alpha^{1/2}\gtrsim N^{-3s+1-3\theta/2}.\]
Since $s<1/3$, we can choose $\theta>0$ such that
$-3s+1-3\theta/2>0$ and it follows that
$\|v\|_{H^s}\rightarrow+\infty$. $\Box$

% The Appendices part is started with the command \appendix;
% appendix sections are then done as normal sections
% \appendix

% \section{}
% \label{}

% Bibliographic references with the natbib package:
% Parenthetical: \citep{Bai92} produces (Bailyn 1992).
% Textual: \citet{Bai95} produces Bailyn et al. (1995).
% An affix and part of a reference:
%   \citep[e.g.][Ch. 2]{Bar76}
%   produces (e.g. Barnes et al. 1976, Ch. 2).

\bibliographystyle{plain}
\bibliography{gbo}

%\begin{thebibliography}{}

% \bibitem[Names(Year)]{label} or \bibitem[Names(Year)Long names]{label}.
% (\harvarditem{Name}{Year}{label} is also supported.)
% Text of bibliographic item

%\bibitem[]{}

%\end{thebibliography}

\end{document}